
\font\bigbold=cmbx12

\font\smallheader=cmssbx10 
\font\bigheader=cmssbx10


\font\tensc=cmcsc10

\font\eightpt=cmr8
\font\ninept=cmr9


\font\ninett=cmtt9
\font\nineit=cmti9

\def\leftrighttop#1#2{
  \headline{\ifnum\pageno=1\hfil\else{\ninept #1 \hfil #2}\fi}
}

\def\firstnopagenum{
  \footline{\ifnum\pageno=1 \hfil \else \hfil{\rm \number\pageno}\hfil\fi}
}

\def\maketitle#1#2#3#4{
  \centerline {\bigbold #1}
  \medskip
  \centerline {\eightpt #2}
  \medskip
  \centerline {\tensc #3}
  \medskip
  \centerline {\tensc #4}
  \bigskip
}


\outer\def\floattext#1 #2. #3\par{
  $$
  \vbox{
    \hsize #1 true in
    \noindent{\bf #2.}\enskip #3
  }
  $$
}


\def\section#1\par{
  \bigskip\vskip\parskip
  \leftline{\bigheader#1}\nobreak\medskip\noindent
}
\def\subsection#1\par{
  \bigskip\vskip\parskip
  \leftline{\smallheader#1}\nobreak\medskip\noindent
}

\def\csection#1\par{
  \bigskip\vskip\parskip
  \centerline{\smallheader#1}\nobreak\medskip\noindent
}

\def\rsection#1\par{
  \bigskip\vskip\parskip
  \rightline{\smallheader#1}\nobreak\medskip\noindent
}

\newskip\ttglue
\ttglue=.5em plus.25em minus.15em

{\obeylines\gdef\startdisplay#1
  {\catcode`\^^M=5$$#1\halign\bgroup\indent##\hfil&&\qquad##\hfil\cr}}
\outer\def\enddisplay{\crcr\egroup$$}

\chardef\other=12
\def\ttverbatim{\begingroup \catcode`\\=\other \catcode`\{=\other
  \catcode`\}=\other \catcode`\$=\other \catcode`\&=\other
  \catcode`\#=\other \catcode`\%=\other \catcode`\|=\other
  \catcode`\_=\other \catcode`\^=\other
  \obeyspaces \obeylines \tt}

{\obeyspaces\gdef {\ }} 

\outer\def\begintt{$$\let\par=\endgraf \ttverbatim \parskip=0pt
  \catcode`\~=0 \rightskip=-5pc \ttfinish}
{\catcode`\~=0 ~catcode`~\=\other 
  ~obeylines 
  ~gdef~ttfinish#1^^M#2\endtt{#1~vbox{#2}~endgroup$$}}



\font\tenfrak=eufm10
\font\sevenfrak=eufm7
\font\fivefrak=eufm5
\newfam\frakfam
\textfont\frakfam=\tenfrak
\scriptfont\frakfam=\sevenfrak
\scriptscriptfont\frakfam=\fivefrak


\def\janksc#1#2 {#1{\eightpt#2}}
\def\jankscsp#1#2 {#1{\eightpt#2}\ }
\def\scproclaim#1.#2\par{\noindent\jankscsp #1.\enspace{\it#2\par}}

\outer\def\articleref #1; #2; #3; #4; #5; #6\par{#1, ``#2,'' {\sl #3}, {\bf #4}, #5, #6.}

\font\ninerm=cmr9
\font\eightrm=cmr8
\font\sixrm=cmr6

\font\ninei=cmmi9 
\font\eighti=cmmi8  
\font\sixi=cmmi6

\font\ninesy=cmsy9 
\font\eightsy=cmsy8
\font\sixsy=cmsy6

\font\nineex=cmex9
\font\eightex=cmex8
\font\sevenex=cmex7

\font\ninebf=cmbx9
\font\eightbf=cmbx8
\font\sixbf=cmbx6

\font\twelvett=cmtt12  \hyphenchar\twelvett=-1  
\font\tensltt=cmsltt10  \hyphenchar\tensltt=-1
\font\ninett=cmtt9  \hyphenchar\ninett=-1
\font\ninesltt=cmsltt10 at 9pt  \hyphenchar\ninesltt=-1
\font\eighttt=cmtt8  \hyphenchar\eighttt=-1
\font\seventt=cmtt8 scaled 875  \hyphenchar\seventt=-1

\font\ninesl=cmsl9
\font\eightsl=cmsl8

\font\nineit=cmti9
\font\eightit=cmti8

\def\ninepoint{\def\rm{\fam0\ninerm}%
  \textfont0=\ninerm \scriptfont0=\sixrm \scriptscriptfont0=\fiverm
  \textfont1=\ninei \scriptfont1=\sixi \scriptscriptfont1=\fivei
  \textfont2=\ninesy \scriptfont2=\sixsy \scriptscriptfont2=\fivesy
  \textfont3=\nineex \scriptfont3=\sevenex \scriptscriptfont3=\sevenex
  \def\it{\fam\itfam\nineit}%
  \textfont\itfam=\nineit
  \def\sl{\fam\slfam\let\ninett=\ninesltt\ninesl}%
  \textfont\slfam=\ninesl
  \def\bf{\fam\bffam\ninebf}%
  \textfont\bffam=\ninebf \scriptfont\bffam=\sixbf
   \scriptscriptfont\bffam=\fivebf
  \def\tt{\fam\ttfam\ninett}%
  \let\sltt=\error
  \textfont\ttfam=\ninett
  \normalbaselineskip=11pt
  \def\bigfences{\textfont3=\tenex}%
  \let\big=\ninebig
  \let\Big=\nineBig
  \let\bigg=\ninebigg
  \let\Bigg=\nineBigg
}

\def\eightpoint{\def\rm{\fam0\eightrm}%
  \textfont0=\eightrm \scriptfont0=\sixrm \scriptscriptfont0=\fiverm
  \textfont1=\eighti \scriptfont1=\sixi \scriptscriptfont1=\fivei
  \textfont2=\eightsy \scriptfont2=\sixsy \scriptscriptfont2=\fivesy
  \textfont3=\eightex \scriptfont3=\sevenex \scriptscriptfont3=\sevenex
  \def\it{\fam\itfam\eightit}%
  \textfont\itfam=\eightit
  \def\sl{\fam\slfam\eightsl}%
  \textfont\slfam=\eightsl
  \def\bf{\fam\bffam\eightbf}%
  \textfont\bffam=\eightbf \scriptfont\bffam=\sixbf
   \scriptscriptfont\bffam=\fivebf
  \def\tt{\fam\ttfam\eighttt}%
  \let\sltt=\error
  \textfont\ttfam=\eighttt
  \def\oldstyle{\fam\@ne\eighti}%
  \normalbaselineskip=9pt
  \def\bigfences{\textfont3=\nineex}%
  \let\big=\eightbig
  \let\Big=\eightBig
  \let\bigg=\eightbigg
  \let\Bigg=\eightBigg
  \normalbaselines\rm}

\def\caption Fig. #1. #2.{{\ninepoint{\bf Fig.\ #1.}\enspace \ninerm#2.}}


\def\xskip{\hskip 7pt plus 3pt minus 4pt}

\def\proof{\medbreak\noindent{\it Proof.}\xskip\ignorespaces}

\def\slug{\quad\hbox{\kern1.5pt\vrule width2.5pt height6pt depth1.5pt\kern1.5pt}\medskip}
\def\noskipslug{\quad\hbox{\kern1.5pt\vrule width2.5pt height6pt depth1.5pt\kern1.5pt}}

\newdimen\algindent
\newif\ifitempar \itempartrue 
\def\algindentset#1{\setbox0\hbox{{\bf #1.\kern.25em}}\algindent=\wd0\relax}
\def\algbegin #1 #2{\algindentset{#21}\alg #1 #2} 
\def\aalgbegin #1 #2{\algindentset{#211}\alg #1 #2} 
\def\alg#1(#2). {\medbreak 
  \noindent{\bf#1}({\it#2\/}).\xskip\ignorespaces}
\def\algstep#1.{\ifitempar\smallskip\noindent\else\itempartrue
  \hskip-\parindent\fi
  \hbox to\algindent{\bf\hfil #1.\kern.25em}%
  \hangindent=\algindent\hangafter=1\ignorespaces}


\def\N{{\bf N}}

\def\R{{\bf R}}

\def\M{{\cal M}} 
 
\def\S{{\cal S}}


\def\ex{\exists}
\def\e{\epsilon}


\def\argmax{\mathop{\rm arg\,max}}

\def\oldno#1{\eqno({\oldstyle#1})}

\outer\def\parenproclaim #1 (#2).#3\par{\medbreak
  \noindent{\bf #1}\enspace\rm({\it #2\/}).\nobreak\ignorespaces{\sl #3\par}
  \ifdim\lastskip<\medskipamount \removelastskip\penalty55\medskip\fi}



\input epsf
\input eplain
\input color
\input soul.sty


\enablehyperlinks
\def\noblueboxes{\special{ps:[/pdfm { /big_fat_array exch def big_fat_array 1 get 0 0 put big_fat_array 1 get 1 0 put big_fat_array 1 get 2 0 put big_fat_array pdfmnew } def}}

\noblueboxes

\immediate\write16{No file \jobname.bbl.} 

\ifpdf 
  \hlopts{bwidth=0} 
  \pdfoutline goto name {sec1} {Introduction}%
  \pdfoutline goto name {sec2} {Automorphisms and Probabilities}%
  \pdfoutline goto name {sec3} {Multiplicities and Correction Factors}%
  \pdfoutline goto name {sec4} {Estimating the Root}%
  \pdfoutline goto name {sec5} {Applications to k-ary and Cayley Trees}%
  \pdfoutline goto name {sec6} count -5 {Tools for Computation on the Conditional Galton--Watson Tree}
    \pdfoutline goto name {conditional} {Events on the conditional tree}%
    \pdfoutline goto name {sums} {Sums of independent random variables}%
    \pdfoutline goto name {numbernodes} {The number of nodes of degree i}%
    \pdfoutline goto name {maximal} {The maximal degree}%
    \pdfoutline goto name {weighted} {Weighted sums}%
    
  \pdfoutline goto name {sec7} count -2 {Probability of Correctness of the Maximum-Likelihood Estimator}
    \pdfoutline goto name {without} {Distributions without special integers}%
    \pdfoutline goto name {with} {Distributions with special integers}%
    
  \pdfoutline goto name {sec8} {Further Examples}
    
\else 
  \special{ps:[/PageMode /UseOutlines /DOCVIEW pdfmark}
  \special{ps:[/Count -0 /Dest (sec1) cvn /Title (Introduction) /OUT pdfmark} 
  \special{ps:[/Count -0 /Dest (sec2) cvn /Title (Automorphisms and Probabilities) /OUT pdfmark}
  \special{ps:[/Count -0 /Dest (sec4) cvn /Title (Estimating the Root) /OUT pdfmark}
  \special{ps:[/Count -0 /Dest (sec5) cvn /Title (Applications to k-ary and Cayley Trees) /OUT pdfmark}
  \special{ps:[/Count -5 /Dest (sec6) cvn /Title (Tools for Computation on the Conditional Galton--Watson Tree) /OUT pdfmark}
     \special{ps:[/Count -0 /Dest (conditional) cvn /Title (Events on the conditional tree) /OUT pdfmark} 
     \special{ps:[/Count -0 /Dest (sums) cvn /Title (Sums of independent random variables) /OUT pdfmark} 
     \special{ps:[/Count -0 /Dest (numbernodes) cvn /Title (The number of nodes of degree i) /OUT pdfmark} 
     \special{ps:[/Count -0 /Dest (maximal) cvn /Title (The maximal degree) /OUT pdfmark} 
     \special{ps:[/Count -0 /Dest (weighted) cvn /Title (Weighted sums) /OUT pdfmark} 

  \special{ps:[/Count -2 /Dest (sec7) cvn /Title (Probability of Correctness of the Maximum-Likelihood Estimator) /OUT pdfmark}
     \special{ps:[/Count -0 /Dest (without) cvn /Title (Distributions without special integers) /OUT pdfmark}
     \special{ps:[/Count -0 /Dest (with) cvn /Title (Distributions with special integers) /OUT pdfmark}

  \special{ps:[/Count -0 /Dest (sec8) cvn /Title (Further Examples) /OUT pdfmark}
\fi

\firstnopagenum  
\def\aut{\mathop{\hbox{\rm Aut}}\nolimits}
\def\stab{\mathop{\hbox{\rm Stab}}\nolimits}
\def\mult{M}
\def\prob{\mathop{\hbox{\rm Prob}}\nolimits}
\def\perm{\mathop{\hbox{\rm Perm}}\nolimits}

\def\Bin{\mathop{\hbox{\rm Binomial}}\nolimits}
\def\Geo{\mathop{\hbox{\rm Geometric}}\nolimits}
\def\Pos{\mathop{\hbox{\rm Poisson}}\nolimits}
\def\Uni{\mathop{\hbox{\rm Uniform}}\nolimits}
\def\boldlabel#1. {\noindent{\bf #1.\enspace}}
\def\pr{\mathop{\hbox{\bf P}}\nolimits}
\def\ex{\mathop{\hbox{\bf E}}\nolimits}
\def\var{\mathop{\hbox{\bf V}}\nolimits}

\def\sc{\ninerm} 

\def\eqref#1{({\oldstyle{#1}})}

\def\comp#1{{#1}^c}
\def\mle{{\cal C}}
\def\one{\mathop{\bf 1}\nolimits}
\def\ind#1{\one_{[#1]}}
\def\bigind#1{\one_{\big[#1\big]}}
\def\eps{\epsilon}
\def\bigmid{\;\Big|\;}
\def\nonameoldno{\oldno{\the\eqcount}\global\advance\eqcount by 1} 

\font\mathbold=cmmib10

\newcount\thmcount  
\thmcount=1
\newcount\sectcount  
\sectcount=1
\newcount\figcount  
\figcount=1
\newcount\eqcount  
\eqcount=1

\definecolor{green}{rgb}{0.1,0.5,0.2}
\definecolor{red}{rgb}{0.5,0.1,0.1}
\definecolor{black}{rgb}{0.0,0.0,0.0}
\def\rev#1{\textcolor{black}{#1}}
\def\red#1{\textcolor{red}{}}



\magnification=\magstephalf
\hoffset=40pt \voffset=28pt
\hsize=29pc  \vsize=45pc  \maxdepth=2.2pt  \parindent=19pt
\nopagenumbers
\def\leftheadline{\folio\hfil{\eightpoint BRANDENBERGER, DEVROYE, AND GOH}\hfil}
\def\rightheadline{\hfil{\eightpoint ROOT ESTIMATION IN GALTON-WATSON TREES}\hfil\folio}
\headline={\ifodd\pageno{\ifnum\pageno<2\hfil\else\rightheadline\fi}\else\leftheadline\fi}

\maketitle{Root estimation in Galton--Watson trees}{}{Anna M. Brandenberger, Luc Devroye, {\rm and} Marcel K. Goh}{{\sl School of Computer Science, McGill University}}
\medskip
\floattext4.6 {\ninebf Abstract}. {\baselineskip=12pt \ninerm Given only the free-tree structure of a tree, the root estimation problem asks if one can guess which of the free tree's nodes is the root of the original tree. We determine the maximum-likelihood estimator for the root of a free tree when the underlying tree is a size-conditioned Galton--Watson tree and calculate its probability of being correct.}
\smallskip
\noindent{\ninebf Keywords.}\enskip\ninerm
Root estimation, Galton--Watson trees, maximum-likelihood methods, probabilistic analysis. 
\baselineskip=14pt

\section\the\sectcount. Introduction 
\hldest{xyz}{}{sec\the\sectcount}

{\tensc Trees are the most} ubiquitous nonlinear structures in computer science. There are two different, equally important, notions of a tree. The first is the {\it unrooted} or {\it free} tree, which is a connected \rev{unlabelled} acyclic graph, and the second is the {\it rooted} tree, in which a single node is distinguished as the root and each edge has a direction from a child to its parent (so all edges point towards the root). Any free tree can be converted into a rooted tree by choosing a root node and setting all of the edge directions accordingly. Likewise, any rooted tree can be seen as a free tree by ``forgetting'' the directions of the edges. The root estimation problem asks for a method that will recover the root of the underlying rooted tree from the free-tree structure.

Given a free tree of size $n$, uniformly chosen from among all $n$-node free trees of a certain family, an easy strategy would be to pick a node uniformly at random; this estimator has a success probability of $1/n$. There are some trees for which this is the optimal estimator, but we will see that in most cases, we will be able to do much better. Of course, it is easy to cook up a family of trees whose structure ensures that the root can be guessed with certainty every time (an obvious example is the the complete binary tree on $2^n-1$ nodes). In many cases we will not be so fortunate, but often there is an estimator that guesses the root with probability asymptotically equal to $c/n$, where $c>1$. We solve the root estimation problem on conditional Galton--Watson trees and exploit the connection between these trees and various families in the uniform tree model to give a general approach to root estimation.



\medskip 


\boldlabel Background.
Root-finding algorithms have been investigated in the literature, mostly for specific classes of trees. The problem was introduced by Haigh~\cite{haigh_1970} in the context of uniform attachment trees, and this work obtains a maximum-likelihood estimate of the root along with the probability of correctness of this estimate as a function of the size of the tree. More recently, Bubeck et.\ al.\ \noblueboxes\cite{bubeck2017finding} show that on uniform attachment and preferential attachment trees, one can construct a confidence set of nodes containing the root, where this set has size independent of the number of nodes in the graph. The earlier work by Shah and Zaman~\cite{shah2011rumors} in network analysis estimates the source of a rumour in a social network under the susceptible-infected-recovered (SIR) model for viral epidemics, which can be viewed as uniform attachment on a background graph. Their estimation is based on the rumour centrality metric, a notion which is explored in further work reviewed in~\cite{shelke2019source}. For instance, Shah et.\ al.\ \cite{shah2016finding} extend their previous result to more generic classes of trees including $d$-regular trees and geometric trees, and further show that their rumour centrality estimator correctly detects the source in Galton--Watson trees with a strictly positive probability. 

In a similar line of work to~\cite{bubeck2017finding}, including some follow up work, authors investigate uniform attachment and preferential attachment trees initialized with an original seed tree ~\cite{bubeck2017trees,bubeck2015influence,curien2015scaling,devroye2019discovery,khim2016confidence,lugosi2019finding}. The authors here seek to determine the original seed of a given graph, and study the influence of this seed and its properties on the structure of the graph as it grows. 
Recent work by Crane et.\ al.\ \cite{crane2020inference} considers {\it shape-exchangable} trees, which encompass the aforementioned models such as uniform attachment, linear preferential attachment, and uniform attachment on a $d$-regular tree, and expand on the ideas of~\cite{bubeck2017finding} and~\cite{khim2016confidence} to provide algorithms for explicitly constructing a confidence set containing the root.

\medskip 
\special{ps:[/pdfm { /big_fat_array exch def big_fat_array 1 get 0 0 put big_fat_array 1 get 1 0 put big_fat_array 1 get 2 0 put big_fat_array pdfmnew } def}
\boldlabel The Galton--Watson model. A Galton--Watson tree~\cite{athreya1972branching} with offspring distribution $\xi$ is a rooted ordered tree in which every node has $i$ children with probability $p_i = \pr\{\xi = i\}$. It is a well-known result that when $\ex\{\xi\} \leq 1$, the tree is finite almost surely, except when $p_1 = 1$ and all other $p_i$ are zero. The Galton--Watson branching process was first studied in 1845 by I.~J.~Bienaym\'e~\cite{bienayme1845}, who was interested in the disappearance of family names, and it derives its name from F.~Galton and H.~W.~Watson~\cite{galtonwatson1874}, who studied the same phenomenon in England in 1874. In their model, nodes correspond to individuals in a population and $p_i$ is the probability that an individual passes the family name down to $i$ children. If the process results in a finite Galton--Watson tree, this means the family name goes extinct after some number of generations. We will consider finite, ``critical'' Galton--Watson trees. These are trees for which $\ex\{\xi\} = 1$ and $\var\{\xi\} \in (0,\infty)$; ensuring a nonzero variance rules out the degenerate case $p_1 = 1$.

The Galton--Watson trees that we shall study are conditioned on $|T| = n$, where $|T|$ is the number of nodes in the tree. Conditional Galton--Watson trees were first studied by D. P. Kennedy~\cite{kennedy1975} and a key correspondence was found between offspring distributions of conditional Galton--Watson trees and certain families of ``simply-generated trees''~\cite{meirmoon1978}:
\medskip
\item{i)} When $\xi\sim \Bin(k, 1/k)$, the conditional Galton--Watson tree is a {\it $k$-ary tree}.
\smallskip
\item{ii)} When $\xi\sim \Pos(1)$, we have a {\it Cayley tree}.
\smallskip
\item{iii)} The distribution $p_0 = p_1 = p_2 = 1/3$ generates a random {\it Motzkin tree}, in which every node has $\leq 2$ children whose order is significant.
\smallskip
\item{iv)} A $\Geo(1/2)$ offspring distribution gives rise to a {uniformly random {\it rooted ordered tree}, also known as a} {\it planted plane tree}.
\medskip
This gives us a way to pick uniformly at random from any such family of trees; we simply generate a conditional Galton--Watson tree, which can be done in linear expected time~\noblueboxes\cite{devroye2012simul}. We will derive a root-estimation strategy for each of the aforementioned families of trees as special cases of our main result.

Our mission can be formalized as follows. Let a conditional Galton--Watson tree with $n$ nodes be given and suppose the directions of the edges are erased, i.e., we are shown only the free-tree structure $F_n$. The goal is to develop a strategy that determines the node with the highest likelihood to have been the root of the original Galton--Watson tree. We would also like to know the probability that we are correct. 
\medskip

\newcount\binaryfour
\binaryfour=\figcount
\advance\figcount by 1
\boldlabel A concrete example. It is instructive to work through a small toy example using a na\"\i ve counting method. Suppose the offspring distribution is
$$p_0 = {1\over 4},\qquad p_1 = {1\over 2},\qquad p_2 = {1\over 4},$$
and all other $p_i = 0$. Conditioning on the number of nodes $n$ generates a binary tree uniformly at random. Fig.~\the\binaryfour\ illustrates the 14 possibilities when $n=4$.
\midinsert
\vbox{
  $$\epsfbox{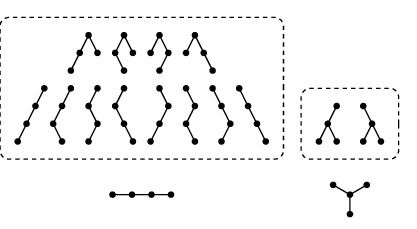}$$
  \vskip -10pt
  \centerline{\caption Fig. \the\binaryfour. The free-tree structure of binary trees with four nodes.}
}
\endinsert
There are only two possible free trees with four nodes and one is much more likely to arise by this process than the other. If {we} are shown a path graph, we are best off choosing one of the endpoints, since an endpoint is the root in 8 of the 12 cases and we will guess the correct endpoint with probability $1/2$ (there are two identical endpoints). In this case, the probability of our guessing correctly is $1/3$. When the free tree is the star graph, we should also choose one of the endpoints, since the central node is never the root. Of course, we can still only be correct with probability $1/3$ because there are three identical endpoints.
\medskip

\boldlabel The probabilistic approach. This family of trees illustrated in Fig.~\the\binaryfour\ was small enough to obtain a maximum-likelihood estimator ({\sc MLE}) by simply counting, but for larger trees and more complex offspring distributions, this will not be feasible. The method we develop will be general and powerful enough to give an {\sc MLE} for the root on conditional Galton--Watson trees with any offspring distribution $p_i$ and any size $n$. We will find that the optimal strategy for picking a root is as follows: 
\medskip
\item{i)} If $p_i > 0$ and $p_{i-1}=0$ for some $i\geq 1$ and there exists a node in the free tree with graph-degree $i$, then only one such node can exist and we select it as our guess. The probability that this node is the root, conditional on its existence {in the free tree}, is 1.
\smallskip
\item{ii)} Otherwise, we choose a node uniformly from the nodes of graph-degree $i$ that maximize $ip_i/p_{i-1}$ (note that there could be multiple integers $i$ for which this ratio is maximal).
\medskip
{Note that computing the {\sc MLE} is computationally easy, and that} the probability of correctness in case (ii) can also be explicitly given. We will also analyze the correctness of the {\sc MLE} as the number of nodes in the tree tends to infinity. 
{Indeed, we show in Theorem~5 that for Galton--Watson trees with offspring distribution satisfying $\sup_{i\geq 1} p_i / p_{i-1} < \infty$ and $0 < \sigma^2 < \infty$, the probability $\pr\{\mle\}$ of the {\sc MLE} being correct satisfies 
$$
\lim_{n\to\infty} n \cdot \pr\{\mle\} = \sup_{i\geq 1} {ip_i \over p_{i-1}}.
$$
Thus, for a large class of tree families for which this supremum is finite, e.g., $k$-ary, Cayley and Motzkin trees, the probability of correctness of the {\sc MLE} decreases linearly with the size of the tree.}

\advance\sectcount by 1
\section\the\sectcount. Automorphisms and Probabilities
\hldest{xyz}{}{sec\the\sectcount}

We start off by establishing some terminology and notation. The setup is as follows. We will denote by $F_n$ a free tree on $n$ nodes. {This is simply an acyclic graph on $n$ vertices, and is {\it a priori} unlabelled, though we may choose labels for the nodes when convenient.} If a node $u$ is selected and the rest of the tree is allowed to hang from it as if by gravity, then we have the {\it $u$-rooted tree}, where the parent of a node is its immediate neighbour in the path towards $u$. 

In the $u$-rooted tree, we define the {\it tree-degree} of a node $v$ to be the number of children of $v$; this is denoted $\deg_u(v)$. The {\it graph-degree} of $v$, written $\deg^*(v)$, is the original degree of $v$ in the free tree $F_n$. For every node $v$ different from $u$ in the $u$-rooted tree, we have $\deg_u(v) = \deg^*(v) - 1$ and $u$ is the only node for which the two degrees are equal. The number of nodes of a given  {tree-}degree $i$ in the $u$-rooted tree is denoted $N_i$; the analogous value for the free tree is denoted $N_i^*$.  {The tree-degree and graph-degree are, in various places, referred to simply as ``degree" (where the context explains which is meant).}

An {\it automorphism} of a free tree $F_n$ is a graph-isomorphism from $F_n$ to itself, i.e., a bijection from $V(F)$ to $V(F)$ that preserves the adjacency structure.  {The group of all such maps is denoted $\aut(F_n)$. We shall define the {\it multiplicity} $M(v)$ of a node $v\in F_n$ to be the size of its orbit under the action of $\aut(F_n)$.}

The notion of free-tree automorphisms is used to define the multiplicity, but in fact the number of automorphisms of a {\it rooted tree} is more pertinent to our problem. Assuming some node $u$ as the root, this is the number of ways that subtrees with the same parent can be permuted amongst themselves while leaving $u$ firmly planted at the top of the tree. In group-theoretic parlance, this is the stabilizer subgroup $\stab(u)$ of the automorphism group of $F_n$.

 {Every Galton-Watson tree is a rooted ordered tree, and we note that if we reorder the children of any given node, we obtain another Galton-Watson tree with exactly the same tree-degree counts, and thus the same probability of occurrence. Repeat this at every node and let $\perm(T)$ be the number of possible such reorderings that one can perform on a given rooted ordered tree $T$; it is clear that there are}
\newcount\autineq
\autineq=\eqcount
\advance\eqcount by 1
$$\prod_{v\in T} \deg_u(v)! \oldno{\the\autineq}$$
{such reorderings.}
But some permutations leave the tree unchanged (if two subtrees of a given node happened to be indistinguishable, then transposing them does not produce a new tree{, in the unordered sense). This happens when, at every node, the reordering only sends children to a slot previously occupied by a node in the same orbit of $\stab(u)$.}
\newcount\correxample
\correxample=\figcount
\advance\figcount by 1
$$\vbox{
  \vskip -25pt
  $$\epsfbox{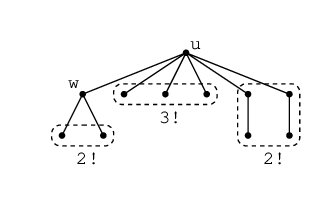}$$
  \vskip -20pt
  \centerline{\caption Fig. \the\correxample. An example tree, in which $\stab(u) = 2!\cdot3!\cdot2! = 96$.}
}$$
For a tree $T$ with root node $u$, we let $\perm(u)$ be the number of distinct unlabelled rooted ordered trees that can be obtained from $T$ by reordering children of nodes.

$$\perm(u) = {1\over \bigl|\stab(u)\bigr|} \prod_v \deg_u(v)\nonameoldno$$

Last but not least, we denote by $\prob(u)$ the Galton--Watson probability of the $u$-rooted tree. Since each node has a probability $p_i$ of having $i$ children, this is given by
$$\prob(u) = \prod_{i=0}^\infty p_i^{N_i}.\nonameoldno$$

Now let $F_n$ be a free tree obtained by removing the parent-child information from a conditional Galton--Watson tree. The probability of a node $u\in F_n$ being the root is the Galton--Watson probability of the $u$-rooted tree times the number {of distinct rooted ordered trees one can obtain via permutations of children. But any node in $u$'s orbit under $\aut(F_n)$ could have been the root of an identical tree, so we must divide by $\mult(u)$. Hence the probability that $u$ is the root is proportional to}
\newcount\probformula
\probformula=\eqcount
\advance\eqcount by 1
$${\prob(u)\perm(u)\over \mult(u)} = {\prob(u)\over \mult(u)\bigl|\stab(u)\bigr|} \prod_v \deg_u(v)! = {\prob(u)\over\bigl|\aut(F_n)\bigr|}\prod_v \deg_u(v)!; \oldno{\the\probformula}$$ 
{one must of course introduce a normalizing factor} to ensure that this is indeed a valid probability distribution. Note that the last equality above is a consequence of the orbit-stabilizer theorem (see, e.g.,~\noblueboxes\cite{hall1959}). Our maximum-likelihood estimator will thus need to choose a node $u$ that maximizes this probability. Given a Galton--Watson offspring distribution, we will denote by $\mle$ the event that the {\sc MLE} is correct for any corresponding free tree of size $n$, and we seek to determine both $\pr\{\mle\}$, the probability of success of the {\sc MLE}, and $\pr\{\mle \mid F_n\}$, the probability of success given a specific free tree $F_n$. Note that
$$\pr\{\mle\} = \ex_{F_n}\bigl\{\pr\{\mle \mid F_n\}\bigr\},\nonameoldno$$
where the expected value is taken over all free trees of size $n$ that could arise by the distribution.

\advance\sectcount by 1
\section\the\sectcount. Estimating the Root
\hldest{xyz}{}{sec\the\sectcount}

We are now ready to prove the first significant result.  {Since $\bigl|\aut(F_n)\bigr|$ does not depend on the choice of root, this boils down to maximizing the quantity $\prob(u)\prod_v \deg_u(v)!$.} The
following theorem shows that this can be done knowing only the offspring distribution and the given free-tree structure $F_n$. To simplify notation, for $i\geq 1$ we define 
$$R_i = {ip_i \over p_{i-1} }.$$
Note that throughout the paper, we will assume that $0/0 = 0$, capturing the cases where both $p_i$ and $p_{i-1}$ are equal to zero.
\newcount\mainthm
\mainthm=\thmcount
\advance\thmcount by 1
\proclaim Theorem \the\mainthm. Given a free tree $F_n$ corresponding to some Galton--Watson tree with offspring distribution $p_i$, the strategy to maximize the probability of picking the original root is to select uniformly from the nodes of graph-degree $i$ that maximize $R_i$, more specifically, defining 
$$ \M =  \max_{j\geq 1} \left \{{R_j} : p_j \neq 0\ \hbox{and there exists}\ u \in F_n\ \hbox{such that}\ \deg^*(u) = j \right \},$$
the maximum-likelihood estimate for picking the root is to choose a node uniformly from the candidate set
$$\Omega = \left \{u \in F_n: \deg^*(u) = i, \; {R_i} = \M \right\}.$$
The probability of success of this maximum-likelihood estimator is 
\medskip
\item{i)} $P\{\mle \mid F_n\} = 1$, if $\M = \infty$;
\smallskip
\item{ii)} when $\M < \infty$, we have
$P\{\mle \mid F_n\} = \M \Big/ \sum_{v \in F_n} R_{\deg^*(v)}.$
\medskip\par

\proof  {The probability that any node $u\in F_n$ is the root is given by the formula~\eqref{\the\probformula}. Thus the goal is to pick a node $u$ that maximizes $\prob(u) \prod_v \deg_u(v)!$.} Suppose we choose some $u$ with $\deg^*(u) = i, \; i \geq 1$.

Note that all the nodes have graph degree one greater than their tree degree, except for the root $u$, where the two degrees are the same. So for all $j \not \in \{i-1, i\}, N_j^u = N_{j+1}^*$ and $N_i = N_{i+1}^* + 1$, $N_{i-1}^u = N_i^* - 1$. We proceed, obtaining
$$\eqalignno{
\prob(u) \prod_v \deg_u(v)! &= \prod_{j=0}^\infty p_j^{N_j^u} \prod_v \deg_u(v)! \cr
&= \prod_{j=0}^\infty p_j^{N_j^u} (j!)^{N_j^u} = \prod_j (j! p_j)^{N_j^u} \cr
&= (i!p_i)^{N_i} \big((i-1)! p_{i-1}\big)^{N_{i-1}^u} \prod_{j\notin \{i, i-1\}} (j! p_j)^{N_j^u} & ({\oldstyle \the\eqcount}) \global\advance\eqcount by 1\cr
&= (i!p_i)^{N_{i+1}^*+1} \big((i-1)! p_{i-1}\big)^{N_{i}^* -1} \prod_{j\notin \{i, i-1\}} (j! p_j)^{N_{j+1}^*} \cr
&= {ip_i\over p_{i-1}} \prod_{ {j=0}}^\infty (j! p_j)^{N_{j+1}^*}. \cr
}$$
The infinite product in the last line is the same for all $u$, so we need only maximize the ratio $R_i$. Considering the constraint that there must be a node of degree $i$ in $F_n$, and the fact that there could be multiple degrees that maximize the required ratio (see the limit of $k$-ary trees as $k\to\infty$ in the following section), there are two cases for the probability of success of this {\sc MLE}. 
\medskip
\item{i)} $\M = \infty$. This case is deceptively simple. If $\M = \infty$, then there exists $i\geq 1$ such that $p_{i-1} = 0$, $p_i \neq 0$, and there is some $u \in F_n$ with $\deg^*(u) = i$. Suppose, towards a contradiction, that this $u$ were not the root. Then there must be some other node $v \neq u$ that is the root, and the $v$-tree degree of $u$ would be $\deg_v(u) = \deg^*(u) - 1 = i-1$. But this is impossible since $p_{i-1} = 0$. So $u$ must be the root. It the only node in the candidate set $\Omega$ and our strategy determines the root correctly with probability $P\{\mle \mid F_n\} = 1$. 
\smallskip
\item{ii)} $\M < \infty$. In this case, since the probability of any node of degree $i$ being the root is proportional to $R_i$, normalizing over all nodes in the free tree $F_n$, we obtain
$$\pr\{\mle\mid F_n\} = \M \Big/ \sum_{v \in F_n} {R_{\deg^*(v)}}.\nonameoldno$$
\medskip
This is exactly the strategy specified in the theorem statement.\slug


\newcount\karysect
\advance\sectcount by 1
\karysect=\sectcount
\section\the\sectcount. Applications to {\mathbold k}{\bf -}ary and Cayley Trees
\hldest{xyz}{}{sec\the\sectcount}

Theorem \the\mainthm\ can be applied to any family of trees that arises as a special case of conditional Galton--Watson trees.  Without any further machinery, we are now able to give an {\sc MLE} for conditional Galton--Watson trees of certain offspring distributions. Recall the computation that we performed on 4-node binary trees in the introduction. We were able to show that the best strategy to guess the root was to choose a random endpoint, which would be successful with probability 1/3. It may come as a surprise that this {\sc MLE} generalises to $k$-ary trees of any size. 

\medskip
\boldlabel Rooted {\mathbold k}-ary trees. In a rooted $k$-ary tree, every node can have up to $k$ children and the placement of the children is important; a node has $k$ ``slots'' in which its children may be placed. As a result, a node can have $i$ children in ${k\choose i}$ ways. When $k = 2$ these trees are often called {\it Catalan trees} because there are ${2n \choose n}/(n+1)$ such trees on $n$ nodes.

We can generate an $n$-node $k$-ary tree uniformly at random by generating a conditional Galton--Watson tree with a $\Bin(k,1/k)$ offspring distribution. Here we have
$$p_i = {k\choose i}\bigg({1\over k}\bigg)^i \bigg({k-1\over k}\bigg)^{k-i}$$
for every $i {\in \{0, \dots, k\}}$, whence
$$R_i = {ip_i\over p_{i-1}} = {i {k\choose i}}{{k\choose i-1}}^{-1} {1\over k}\cdot {k\over k-1} = {k-i+1 \over k-1}.\nonameoldno$$
So, for any free tree $F_n$, the probability of a given node $u$ of degree $\deg^*(u) = i$ being the root is
$$R_i\Big/{\sum_v R_{\deg^*(v)} } = {k-i+1\over \sum_v \big(k - \deg^*(v) + 1\big)} = {k-i+1\over nk - (2n-2) + n} = {k-i+1\over (k-1)n + 2}.\nonameoldno$$
Following the {\sc MLE} strategy, we pick uniformly at random out of the nodes in the free tree with degree $i=1$ (of which at least one is guaranteed to exist). Note that this expression is independent of the shape of the free tree $F_n$, so the probability of success of the {\sc MLE} is
\newcount\karyformula
\karyformula=\eqcount
$$\pr\{\mle\} = \pr\{\mle \mid F_n \} = {k\over (k-1)n + 2}.\nonameoldno$$
\medskip

\boldlabel Cayley trees. From the formula ({\oldstyle{\the\karyformula}}), one can see that for random $k$-ary trees, our advantage decreases as $k$ gets large. Indeed, taking the limit as $k\to \infty$, the $\Bin(k, 1/k)$ distributions approach a $\Pos(1)$ distribution, with $p_i = (e\cdot i!)^{-1}$. This generates the family of {\it Cayley trees}, and in this case,
$${ip_i \over p_{i-1}} = {i\cdot e\cdot (i-1)!\over e\cdot i!} = 1,\nonameoldno$$
so every node is equally likely to be the root. Here there is no better strategy than picking uniformly from all nodes in the tree and the success probability is $1/n$.
\medskip

In both of these cases, $\pr\{\mle \mid F_n\}$ only depends on $n$, and we thus have $\pr\{\mle\} = \pr\{\mle \mid F_n\}$, lending to easy analysis of the {\sc MLE}. This will not be true in all cases, so in the remainder of this paper, we will upgrade the probabilistic technology in our arsenal before reframing the maximum-likelihood estimator and its probability of correctness $\pr\{\mle\}$ for more complex offspring distributions. 

\advance\sectcount by 1 
\section \the\sectcount. Tools for Computation on the Conditional Galton--Watson Tree
\hldest{xyz}{}{sec\the\sectcount}

We would like to be able to analyze the unconditional correctness of the {\sc MLE} $\pr\{\mle\}$ on a Galton--Watson tree with offspring probability $p_i$. {In general, $\pr\{\mle\}$ is a random variable that depends on the free-tree structure of the Galton-Watson tree, and} we will need certain results from the theory of branching processes. This section contains a potpourri of lemmas and small results that will be useful in the upcoming sections and examples. The casual reader may choose to skim through them in anticipation of the main theorems of the next section, returning to enjoy the proofs after seeing the lemmas used in action.
\medskip


\boldlabel Events on the conditional tree. 
\hldest{xyz}{}{conditional}
Let $B$ be some event concerning an unconditional Galton--Watson tree $T$ with offspring distribution $\xi$. We would like to establish useful tools for working with
$$\pr\{B\mid |T| = n\},$$
using the random walk representation of conditional Galton--Watson trees. First, suppose that we number the nodes in $T$ (in {depth-first} preorder, say). Each node $i$ has degree $\xi_i$ and if $\xi_1, \xi_2,\ldots$ are independent and all distributed as $\xi$, then we have
\special{ps:[/pdfm { /big_fat_array exch def big_fat_array 1 get 0 0 put big_fat_array 1 get 1 0 put big_fat_array 1 get 2 0 put big_fat_array pdfmnew } def}
$$\eqalign{
|T| &= \min \{ t > 0 : 1 + (\xi_1 - 1) + \cdots + (\xi_t - 1) = 0\} \cr
&= \min \Big\{ t > 0 : \sum_{i=1}^t \xi_i = t-1\Big\}. \cr
}\nonameoldno$$
Defining two events
$$ A^* = \Big\{ 1 + \sum_{i=1}^t (\xi_i - 1) > 0\ \hbox{for all}\ t < n , \sum_{i=1}^n \xi_i = n-1 \Big\} $$
and
\newcount\eventAeq
\eventAeq=\eqcount
\advance\eqcount by 1
$$ A = \Big\{\sum_{i=1}^n \xi_i = n-1 \Big\}, \oldno{\the\eventAeq}$$
we have, by {Dwass's} cycle lemma~\cite{dwass1969},
$$\pr\{ |T| = n \} = \pr \{A^*\} = {1\over n} \pr\{A\}.\nonameoldno$$
Now, $B$ is an event on $T$, and is thus determined by $\xi_1, \ldots, \xi_n$. If we assume rotation invariance ($B$ remains true if applied to $\xi_i, \xi_{i+1}, \ldots, \xi_n, \xi_{n+1}, \ldots, \xi_{i-1}$ for all $i$), then we obtain, by another use of the cycle lemma,
\newcount\treetoeventeq
\treetoeventeq=\eqcount
\advance\eqcount by 1
$$\eqalign{
\pr\{ B \mid |T| = n\} &= {\pr\{ B \cap |T| = n\}\over\pr\{ |T| = n\}} = {\pr\{B \cap A^*\}\over \pr\{A^*\}} = {\pr\{ B \cap A\} / n\over \pr\{A\}/n}\cr
&= {\pr \{B\cap A\} \over \pr\{A\}} = \pr\{B\mid A\}.
}\oldno{\the\treetoeventeq}$$
This matters because one can study $B$ by simply looking at sequences of i.i.d.\ random variables and without having to worry about trees.
\medskip

\boldlabel Sums of independent random variables.
\hldest{xyz}{}{sums}
We will need two lemmas regarding the sums of random variables; these are well-known and will be given without proof. Let the {\it period} of a random variable $\xi$ be the greatest common divisor of all the $i$'s for which $\pr\{\xi = i\} > 0$. The first of these lemmas is due to B. A. Rogozin~\noblueboxes\cite{rogozin1961} and the statement as well as its proof can be found in~\cite{petrov1975}.

\parenproclaim Lemma A (Rogozin{\rm, 1961}). If $X_1, \ldots, X_n$ are i.i.d.\ random variables and
$$p = \sup_x \pr\{X_1 = x\},$$
then
$$\sup_x \pr\{ X_1 + \cdots + X_n = x\} \leq {\alpha\over \sqrt{n(1-p)}}\nonameoldno$$
for a universal constant $\alpha$.\slug

The following lemma regards sums of i.i.d.\ random variables (e.g.,  as present in our event $A$) and is due to V. F. Kolchin~\cite{kolchin1986}.

\parenproclaim Lemma B (Kolchin{\rm, 1986}). Let $\xi_1,\ldots,\xi_n$ be i.i.d.\ random variables on $[0,\infty)$ of mean 1 and variance $\sigma^2 > 0$. Let the period of $\xi_1$ be
$$ h = \gcd \{ i \geq 1 : p_i > 0\}$$
and let $X$ be the set of all integers $x$ such that $(n+x) \bmod h = 0$. Then
$$\sup_{x\in X} \sqrt n\;\bigg| \pr \{\xi_1 + \cdots + \xi_n = n + x\} - {h\over \sigma\sqrt{2\pi}} e^{-x^2 / 2n\sigma^2}\bigg| \to 0\nonameoldno$$
as $n\to \infty$. If $(n + x)\bmod h \neq 0$, then $\pr\{ \xi_1 +\cdots + \xi_n {= n + x}\} = 0$.\slug
\medskip

\boldlabel The number of nodes of degree {\mathbold i}.
\hldest{xyz}{}{numbernodes}
Recall that we write $N_i$ to indicate the number of nodes of tree-degree $i$ in a Galton--Watson tree. We will show a {result} that as $n$ gets large, the proportion of nodes in the tree of degree $i$ approaches $p_i$. The following lemma is due to Aldous~\cite{aldous1991} and Janson~\cite{janson2016}.

\parenproclaim Lemma C (Aldous{\rm, 1991;} Janson{\rm, 2016}). Let $T_n$ be a conditional Galton--Watson tree with offspring distribution $\xi$ satisfying ${0 <} \sigma^2 < \infty$, and let
$$N_i = \sum_{k=1}^n \ind{\xi_k = i}$$
be the number of nodes of degree $i$ in $T_n$. For any $i$, $N_i/ n \to p_i$ in probability as $n\to \infty$.\slug

 \proof Let $\eps > 0$ be given. Let $A$ be the event that $\sum_{i=1}^n \xi_i = n-1$ and let $B$ be the event that $|N_1 / n = p_i| > \eps$. Note that $B$ is rotation invariant. So we have, by ({\oldstyle{\the\treetoeventeq}}),
 $$\pr\{ B \mid |T| = n\} = \pr\{B\mid A\} = {\pr \{B\cap A\}\over \pr\{A\}} \leq {\pr\{B\}\over \pr\{A\}}.\nonameoldno$$
 Now, by Lemma B,
 $$\pr\{A\} = \pr\Big\{\sum_{i=1}^n \xi_i = n-1\Big\} = {he^{-1/2\sigma n^2} + o(1)\over \sigma\sqrt{2\pi n}} \sim {h\over \sigma\sqrt{2\pi n}},\nonameoldno$$
 where $h$ is the period of $\xi_1$. Also, since $\ex\{N_i\} = np_i$ and
 $$\var\{N_i / n\} = {\var\{\ind{\xi_1 = i}\}\over n} = {p_i(1-p_i)\over n},$$
 we have, by Chebyshev's inequality,
 $$\pr\{B\} \leq {\var \{N_i/n\} \over \eps^2} \leq {p_i(1-p_i)\over n\eps^2},$$
 whence
 $$\pr\{B \mid |T| = n\} \leq {1\over \sqrt n}\bigg({p_i(1-p_i) \sigma \sqrt{2\pi} \over h\eps^2} + o(1)\bigg),\nonameoldno$$
 and the right hand side goes to 0 as $n\to \infty$.\slug
 \medskip

\boldlabel The maximal degree. 
\hldest{xyz}{}{maximal}
Another important random variable is the maximal degree $M_n$ of $T_n$. Because this is rotation-invariant, one can study $M_n$ just as one studies the maximum of independent random variables.

\newcount\maxdegree
\maxdegree=\thmcount
\advance\thmcount by 1
\proclaim Lemma \the\maxdegree. Let $T$ be a conditional  Galton--Watson tree of size $n$ with offspring distribution $\xi$ {whose variance $\sigma^2$ satisfies} $0 < \sigma^2 < \infty$ and let
$$M_n = \max_{1\leq i\leq n} \xi_i$$
be the maximal degree among all the nodes in $T$.
{Fix an integer $x$. Letting $o(1)$ stand for any quantity that tends to $0$ as $n\to\infty$ independent of $x$}, we have
$$\pr\{{M_n\geq x \mid |T| = n }\} \leq \big(1 + o(1)\big)n\pr\{\xi\geq x\}\nonameoldno$$
and
$$\pr\{{M_n\leq x \mid |T| = n }\} \leq \big(\beta + o(1)\big) \exp\big({-n} \pr\{\xi > x\}\big),\nonameoldno$$
for a universal constant $\beta$.

Note that if we have a sequence of $n$ i.i.d.\ random variables $\xi_i$, the same bounds can be derived, without the $\big(1 + o(1)\big)$ and $\big(\beta + o(1)\big)$ factors. This lemma shows that asymptotically, nothing is lost by conditioning on the size of a Galton--Watson tree.

\proof 
{Let $A$ be the event that $\sum_{i=1}^n \xi_i = n-1$.}
We begin by expanding and applying the union bound, with $A$ being the event as in ({\oldstyle{\the\eventAeq}}), obtaining
$$\eqalign{
\pr\{ M_n \geq x \mid |T| = n\} &= \pr\{ M_n \geq x, A\}/ \pr\{A\} \cr
&\leq n\pr \Big\{ \xi_i \geq x,\; \sum_{i=1}^n \xi_i = n-1\Big\} / \pr\{A\} \cr
&= n \sum_{j=x}^\infty \bigg(\pr\Big\{ \xi_1 = j,\; \sum_{i=2}^n \xi_i = n-1-j\Big\} /\pr\{A\}\bigg).\cr
}\nonameoldno$$
Let $h$ be the period of $\xi_1$. By Lemma B, we can proceed as follows:
$$\eqalign{
\pr\{ M_n \geq x \mid |T| = n\} &\leq n\sum_{j=x}^\infty \bigg(p_j \pr\Big\{\sum_{i=2}^n \xi_i = n-1-j\Big\} \bigg/ {h\big(1 + o(1)\big)\over \sigma\sqrt{2\pi n}}\bigg) \cr
&= n\sum_{j=x}^\infty \bigg(p_j {h e^{-j^2/2\sigma^2(n-1)} + o(1)\over \sigma\sqrt{2\pi(n-1)}} \bigg/ {h\big(1 + o(1)\big)\over \sigma\sqrt{2\pi n}}\bigg) \cr
&\leq n\sqrt{n\over n-1} \Big( \sum_{j\geq x} p_j \Big)\big(1 + o(1)\big)
\cr
&\leq \big(1 + o(1)\big)n\pr\{\xi\geq x\}. \cr
}\nonameoldno$$

Next we tackle the lower bound, by an independence argument. First, denoting by $A$ the event that $\sum_{i=1}^n \xi_i = n-1$ as in ({\oldstyle{\the\eventAeq}}), we expand
$$\eqalign{
\pr\{M_n \leq x \mid |T| = n\} &= {\pr\{M_n \leq x, A\} \over \pr\{A\}} \cr
&= \pr\{M_n\leq x\}{\pr \{A \mid M_n\leq x\} \over \pr\{A\}}.\cr
}$$
Well, $\pr\{A \mid M_n\leq x\} = \pr\{ \sum_{i=1}^n \xi_i^* = n-1\}$, where $\xi_1^*,\ldots,\xi_n^*$ are i.i.d.\ with
$$\pr \{\xi_1^* = i\} = \cases {\pr\{\xi_1 = i\}/ \pr\{\xi_1 \leq x\}, & if $i\leq x$; \cr 0, & if $i > x$. \cr}$$
Let $\eta = \min\{i > 0 : p_i > 0\}$. Then, for $x \geq \eta$, we have $p := \max_{i\leq x} p_i / (p_0 + \cdots + p_x) < 1$. Therefore, by Lemma A,
$$\pr\Big\{\sum_{i=1}^n \xi_i^* = n-1\Big\} \leq {\alpha \over \sqrt{n(1-p)}},\nonameoldno$$
for a general constant $\alpha$. Putting
$$\beta = {\alpha \over \sqrt{1-p}}\cdot{\sigma\sqrt{2\pi} \over h},$$
we have, for {$x \geq \eta$},
$$\eqalign{
\pr\{M_n \leq x\mid |T| = n\} &\leq \pr\{M_n\leq x\}\cdot \beta\big(1 + o(1)\big) \cr
&\sim \beta\big(\pr\{\xi \leq x\}\big)^n \cr
&\leq \beta \exp\big({-n}\pr\{\xi\geq x\}\big), \cr
}\nonameoldno$$
On the other hand, if $x < \eta$, then $\pr\{M_n \leq x \mid |T| = n\} = 0$ since $M_n < \eta$ implies $M_n = 0$ and thus $|T| = 1$, which is impossible for $n > 1$. The above bound therefore still holds.
\slug

\boldlabel Weighted sums. 
\hldest{xyz}{}{weighted}
In the derivation of $\pr\{\mle\}$, one encounters the sum
$$\sum_v R_{\deg^*(v)} $$
for a given free tree in the denominator. When the nodes of a conditional Galton--Watson tree are numbered from 1 to $n$ in preorder and each node $i$ produces a number of offspring distributed as $\xi_i$, this sum {is very close to} the random variable
\newcount\Wneqdef
\Wneqdef=\eqcount
$$W_n = \sum_{i=1}^n {1\over p_{\xi_i}} (\xi_i + 1) p_{\xi_i+1} {\ind{p_{\xi_i} \neq 0}}.$$
We give two lemmas that allow us to work with these weighted sums.

\newcount\weightedsum
\weightedsum=\thmcount
\advance\thmcount by 1
\proclaim Lemma~\the\weightedsum. 
Consider the random variable $W_n$ {with} $\xi$ satisfying $0 < \sigma^2 <\infty$ {and $\sup_{i\geq 1, i \not\in \S} p_{i}/p_{i-1} < \infty$, where $\S = \{ i \in \N : p_i > 0, p_{i-1} = 0\}$}{. Defining $\gamma = \sum_{j \not\in \S} jp_j \leq 1$,} we have $W_n/{(\gamma n)} \to 1$ in probability as $n\to \infty$. {Note that if $\S = \emptyset$, then $W_n / n \to 1$.}

\proof Note that
\newcount\meanXeq
\meanXeq=\eqcount
$$\ex\bigg\{ {(\xi + 1)p_{\xi + 1} \over p_\xi} {\ind{p_{j} \neq 0}}\bigg\} = \sum_{j=0}^\infty {p_j\over p_j} (j+1)p_{j+1} {\ind{p_{j} \neq 0}} = \sum_{j=1}^\infty jp_j {\ind{p_{j-1} \neq 0}} = {\gamma},\nonameoldno$$
so $\ex\{W_n\} = {\gamma}n$. {Also,} 
\newcount\variXeq 
\variXeq=\eqcount
$$\eqalign{
\ex\bigg\{\bigg({(\xi + 1)p_{\xi + 1} \over p_\xi} {\ind{p_{j} \neq 0}}\bigg)^2\bigg\} 
&= \sum_{j=0}^\infty {p_j\over {p_j}^2}(j+1)^2 p_{j+1}^2 {\ind{p_{j} \neq 0}} \cr
&\leq {\sup_{j\geq 1, j \not\in \S}} {p_{j}\over p_{j-1}} \sum_{j=0}^\infty (j+1)^2p_{j+1} \cr
&= {\sup_{j\not\in \S}{ p_{j}\over p_{j-1}} (\sigma^2 + {\gamma^2})}.\cr
}\nonameoldno$$
By Chebyshev's inequality, for any arbitrary $\e > 0$
$$\pr\bigg\{ \bigg| {W_n\over {\gamma n} } - 1\bigg| > \eps\bigg\} \leq {\var\{W_n\}\over n^2\eps^2} \leq {(\sigma^2 + {\gamma^2}) {\sup_{j\not\in\S} p_{j}/p_{j-1}} \over n\eps^2}.$$
Therefore, arguing as before and letting $A$ be the event as in ({\oldstyle{\the\eventAeq}}),
$$\pr\bigg\{\bigg|{W_n\over {\gamma n} } - 1\bigg| > \eps \bigmid |T|=n\bigg\} \leq {\pr\big\{ |W_n/ {\gamma n} - 1| > \eps\big\} \over \pr\{A\}} = O\bigg({1\over \sqrt n}\bigg). \noskipslug$$

We would now like to show that $\ex\{{\gamma n}/W_n \mid |T|=n\} \to 1$. This does not follow directly from Lemma~\the\weightedsum, but we shall squeeze it out by means of some well-known inequalities and a little elbow grease.

\newcount\weightedcontd
\weightedcontd=\thmcount
\advance\thmcount by 1
\proclaim Lemma \the\weightedcontd. Under the {same} assumptions {as the previous lemma,} we have
$$\ex\bigg\{ {{\gamma n}\over W_n} \bigmid |T| = n\bigg\} \to 1 \qquad {\hbox{and} \qquad \ex\bigg\{ {{(\gamma n)}^2 \over {W_n}^2} \bigmid |T| = n \bigg\} \to 1} \nonameoldno$$
as $n\to \infty$.

\proof Let $\eps > 0$, and as before let $A$ be the event as in ({\oldstyle{\the\eventAeq}}). First, we observe that
$$\eqalign{
\ex\bigg\{ {{\gamma n}\over W_n} \bigmid |T| = n\bigg\} &\geq {{\gamma n}\over {\gamma n}(1+\eps)} \cdot {\pr\{W_n < {\gamma n}(1+\eps),\; A\}\over \pr\{A\}} \cr
&= {1\over 1+\eps} \left( 1 - {\pr\{W_n \geq {\gamma n}(1+\eps),\; A\} \over \pr\{A\}}\right) \cr
&\geq {1\over 1+\eps} - O\bigg({1\over \sqrt n}\bigg), \cr
}\nonameoldno$$
since $W_n/{\gamma n}\to 1$ in probability and $\pr\big\{ W_n \geq {\gamma n}(1+\eps) \big\} = O(1/n)$, by the previous lemma. {Similarly we have
$$\eqalign{
\ex\left\{ {{(\gamma n)}^2 \over {W_n}^2} \bigmid |T| = n\right\} &\geq {{(\gamma n)}^2 \over {(\gamma n)}^2(1+\e)^2} {\pr\{W_n \geq {\gamma n}(1+\e), A\} \over \pr\{A\}}\cr
&\geq {1 \over (1+\e)^2} - O\left({1 \over \sqrt{n}}\right).
}$$}
It remains to show that $\ex\big\{ {\gamma n} / W_n \;\big|\; |T| = n\big\} \leq 1 + o(1)$ and similarly for $(\gamma n)^2 / {W_n}^2$. To that end, note that
$$W_n \geq \sum_{i=1}^n \ind{\xi_i = 0} \cdot {1\over p_0} p_1.\nonameoldno$$
Letting $N_0 = \sum_{i=1}^n \ind{\xi_i = 0}$, we remark that $N_0\sim \Bin(n, p_0)$ and apply Hoeffding's bound to obtain, for $\delta < \min\{p_0, 1-p_0\}$,
$$\pr\big\{ |N_0 - np_0| > \delta n\big\} \leq 2e^{-2n\delta^2}.$$
{We choose} $\delta = \eps/n^{1/\eps}$. {Then, by rotation-invariance of $W_n$, we have}
$$\ex\bigg\{ {{\gamma n}\over W_n} \bigmid |T| = n\bigg\} = {\ex\big\{ ({\gamma n}/W_n)\one_A \big\} \over \pr\{A\}}\nonameoldno$$
and
$$\ex\bigg\{ {(\gamma n)^2/{W_n}^2} \bigmid |T| = n\bigg\} = {\ex\big\{ ((\gamma n)^2/{W_n}^2)\one_A \big\} \over \pr\{A\}},\nonameoldno$$
{where} $\pr\{A\} = \Theta(1/\sqrt n)$. {Recall that $a_n = \Theta(b_n)$ denotes the existence of constant real numbers $c, d > 0$ such that for large enough $n$, $a_n \leq c b_n$ and $a_n \geq d b_n$.} Also,
\newcount\expectedtwoterms
\expectedtwoterms=\eqcount
$$\eqalign{
\ex\bigg\{ {{\gamma n}\over W_n} \one_A \bigg\} &\leq \ex\bigg\{{{\gamma n}\over (1-\eps){\gamma n}} \one_A \bigg\} + \ex\bigg\{ {{\gamma n}\over W_n} \ind{W_n\leq (1-\eps){\gamma n}} \one_A \bigg\} \cr
&\leq {1\over 1-\eps}\pr\{A\} + \ex\bigg\{ {p_0\over p_1} \cdot {{\gamma n}\over N_0} \cdot \ind{N_0 \leq np_0/2} \cdot \one_A\bigg\}\cr
&\qquad\qquad+ \ex \bigg\{{p_0\over p_1}\cdot {2{\gamma n}\over np_0} \cdot \ind{W_n\leq (1-\eps) {\gamma n}}\bigg\}.\cr
}\nonameoldno$$
Letting $E_1$ and $E_2$ denote the two expectation terms on the right-hand side, we note that since $A$ implies that $N_0\geq 1$,
$$E_1 \leq {p_0 \over p_1} {\gamma n} \pr \{N_0\leq np_0/2\} \leq {p_0\over p_1} 2n\exp\big({-2}{\gamma n}(p_0/2)^2\big).$$
Furthermore,
$$E_2 = {2 {\gamma}\over p_1} \pr\big\{W_n \leq (1-\eps){\gamma n}\big\} = O\bigg({1\over n}\bigg)$$
follows from Chebyshev's inequality, just as in the proof of Lemma~\the\weightedsum. This implies that
$$\ex\bigg\{ {{\gamma n}\over W_n} \bigmid |T| = n\bigg\} \leq {1\over 1-\eps} + {O(1/n) \over O(1/\sqrt n)} = {1+o(1)\over 1-\eps},\nonameoldno$$
and we are done since $\eps$ was chosen arbitrarily. {For the $(\gamma n)^2 / {W_n}^2$ case, we proceed the same way to obtain
$$\eqalign{ 
\ex\left\{ {{(\gamma n)}^2 \over {W_n}^2}\one_A \right\} 
&\leq {1 \over (1-\e)^2} \pr\{A\} + \ex\bigg\{ {(\gamma n)^2 \over {W_n}^2 } \ind{W_n\leq (1-\eps)(\gamma n)^2 } \one_A \bigg\} \cr 
&\leq {1 \over (1-\e)^2} \pr\{A\} + \ex\bigg\{ {p_0^2\over p_1^2} \cdot {(\gamma n)^2\over {N_0}^2} \cdot \ind{N_0 \leq np_0/2} \cdot \one_A\bigg\}\cr
&\qquad\qquad+ \ex \bigg\{{p_0^2\over p_1^2}\cdot {4 (\gamma n)^2\over n^2p_0^2} \cdot \ind{W_n\leq (1-\eps) {\gamma n}}\bigg\}\cr 
&\leq {1 \over (1-\e)^2} \pr\{A\} + {p_0^2 \over p_1^2} (\gamma n)^2 \pr \{N_0\leq np_0/2\}\cr
&\qquad\qquad+ {4 \gamma^2 \over p_1^2} \pr\{W_n \leq (1-\e) \gamma n\},\cr 
}$$
completing the proof in the same manner.}
\slug


\advance\sectcount by 1 
\section\the\sectcount.  Probability of Correctness of the Maximum-Likelihood Estimator
\hldest{xyz}{}{sec\the\sectcount}

We begin by setting up a few definitions to better deal with the two cases mentioned in Theorem~\the\mainthm\ in the large $n$ limit. Using this notation, we reformulate our maximum-likelihood estimator for the root, and compute its expected probability of correctness $\pr\{\mle\}$.

Let an offspring distribution be fixed. If $p_i > 0$ and $p_{i-1} = 0$ for some positive integer $i$, we say that $i$ is a {\it special integer} and we call a node in the free tree with graph degree $i$ a {\it special node}. Remember that finding a special node is akin to hitting the jackpot for the {\sc MLE}. If $i$ is a special integer and some node $v$ in a free tree has graph-degree $i$, then $v$ is the root with probability 1. We denote the set of all special integers by $\S$. Note that $i=1$ is never special, since $p_0 > 0$. We group all non-special integers $i$ into equivalence classes $\{J_k\}_{k\geq 1}$ according to the equivalence
$$i\sim j\quad\hbox{if and only if}\quad{ip_i\over p_{i-1}} = {jp_j\over p_{j-1}}.$$
As before, we let $R_i = ip_i/p_{i-1}$ but for convenience, we will allow the notation $R_{J_k}$, which equals $R_i$ for any $i\in J_k$. 
%
%
Lastly, we let $N_{J_k}$ denote the number of nodes in the tree whose graph-degree belongs in the equivalence class $J_k$; recalling that $N_i^*$ is the number
of nodes with graph-degree $i$, we have
$$N_{J_k} = \sum_{i\in J_k} N^*_i.$$

\medskip 

\boldlabel The maximum-likelihood estimator. With these new definitions, we can formally redescribe the {\sc MLE} and the probability of correctness. Given a free tree $F_n$ of size $n$ corresponding to a Galton--Watson tree with offspring distribution $p_i$, we guess the root as follows:
\medskip 
\item{i)} Let $S_n$ denote the event that there exists a special node in a given free tree $F_n$. If $S_n$ occurs, then select this special node. In this case, 
$$ \pr\{\mle \mid F_n\} \one_{S_n} = \one_{S_n}.$$
\smallskip
\item{ii)} Otherwise, {let $\comp S_n$ denote the complement of $S_n$ which occurs} if there are either no special integers in the distribution or no nodes with the corresponding degrees in the free tree. {On this event,} select a node uniformly at random from the class $J_\lambda$, where $$ \lambda = \argmax_{k\notin\S}\{R_{J_k}: N_{J_k} > 0\},$$
where we note that this maximum can be taken since there are at most $n$ non-empty equivalence classes. In this case,
$$\pr\{\mle \mid F_n\} {\one_{\comp S_n}} = {R_\lambda \over \sum_{k} N_{J_k} R_{J_k}} {\one_{\comp S_n}}.\nonameoldno$$
\medskip

\special{ps:[/pdfm { /big_fat_array exch def big_fat_array 1 get 0 0 put big_fat_array 1 get 1 0 put big_fat_array 1 get 2 0 put big_fat_array pdfmnew } def}
\boldlabel Distributions without special integers.
\hldest{xyz}{}{without}
We first consider the well-behaved (and more common) case in which there exist no special integers in the Galton--Watson distribution $p_i$. The following theorem will require the notion of Kesten's limit tree~\noblueboxes\cite{kesten1986}, which we will briefly describe. Recall that we are working with an offspring distribution $\xi$ for which $\ex\{\xi\} = \sum_{i\geq 1} ip_i = 1$. So if $\zeta$ is the random variable with $\pr\{\zeta = i\} = ip_i$ for all $i\geq 1$, then $\zeta$ is a valid offspring distribution as well. Kesten's limit tree $T_\infty$ is an infinite tree consisting of a central spine of nodes, one on each level, that each produce $\zeta$ children. Nodes that are not on the spine are the root of an unconditional Galton--Watson tree with distribution $\xi$ (each of these is finite with probability 1). Let $\tau(T, h)$ denote the tree $T$, limited to levels $0,\ldots, h$. Kesten's limit tree is important to us because for all $h$ and all infinite trees $t$, {a Galton--Watson tree $T_n$ conditioned to be of size $n$ converges locally to it in the following sense:}
$$\lim_{n\to\infty} \pr\big\{\tau(T_n, h) = \tau(t,h)\big\} = \pr\big\{\tau(T_\infty, h) = \tau(t,h)\big\}.\nonameoldno$$

\newcount\mainmainthm
\mainmainthm=\thmcount
\advance\thmcount by 1
\proclaim Theorem \the\mainmainthm. Given a random free tree of size $n$ corresponding to a Galton--Watson tree with offspring distribution $p_i$ with $0 < \sigma^2 < \infty$ and $\sup_{i\geq 1} p_i/p_{i-1} < \infty$. Then the probability of the {\sc MLE} being correct satisfies 
$$\lim_{n\to\infty} n\cdot \pr\{\mle\} = \sup_{i\geq 1} {ip_i \over p_{i-1}}.\nonameoldno$$
Note that this could be infinity.

\proof 
Let $\lambda$ indicate the equivalence class chosen by the {\sc MLE}, as described above. First, we prove the upper bound:
$$\pr\{\mle\} = \ex\big\{\pr\{\mle \mid F_n\}\big\} = \ex\left\{{R_\lambda \over \sum_k N_{J_k} R_{J_k}}\right\} \leq \sup_{i\geq 1} R_i \ex\left\{{1\over \sum_k N_{J_k} R_{J_k}} \right\},$$
where we note that $\sum_k N_{J_k} R_{J_k} = \sum_v R_{\deg^*(v)}$ corresponds, up to a $O(1)$ error, to the random variable $W_n$ from Lemma~\the\weightedcontd, which gave us that $\ex\{n/W_n \mid |T| = n\} \to 1$ {as $\gamma = 1$}. We can thus conclude that 
$$\limsup_{n\to\infty} n\pr\{\mle\} \leq \sup_{i\geq 1} R_i. \nonameoldno$$
Before moving to the lower bound, let us first show that for any degree $i \geq 1$ such that $p_i > 0$, as $n \to \infty$,
$$\pr\{N_i^* = 0\} \to 0.$$ 

Note that by Lemma~C, for any conditional Galton--Watson tree corresponding to the free tree of size $n$ rooted at a node $u$, for all $i$, $N_i / np_i \to 1$ in probability. Furthermore, since we assumed that our distribution has no special integers, for any degree $i$ such that $p_i > 0$, we also have $p_{i-1} > 0$. This yields, for any $i \geq 1$, 
$$\eqalign{
\pr\{N_i^* = 0\} &= \pr\{N_i^* = 0 , \deg^*(u) \not \in  \{i, i-1\} \}\cr
&\qquad\qquad+ \pr\{N_i^* = 0 , \deg^*(u) = i-1 \} + \pr\{N_i^* = 0 , \deg^*(u) = i\}  \cr
&= \pr\{N_{i-1}^u = 0 \mid \deg^*(u) \not \in \{i, i-1\}\}\pr\{\deg^*(u) \not \in \{i, i-1\}\} \cr 
&\qquad\qquad + \pr\{N_{i-1}^u - 1 = 0 \mid \deg^*(u) = i-1\}\pr\{\deg^*(u) = i-1\}, \cr 
}$$
which goes to 0. This follows from the fact that, as $n$ gets large and the conditional Galton--Watson tree converges {locally} to Kesten's limit tree, $\pr\{\deg^*(u) = i\} = ip_i + o(1)$. Note that in the above argument, the random variables $N^*_i$, $N^u_i$ and $\deg^*(u)$ all depend on $n$, but we avoid double-indexing for clarity of notation.

Now for the lower bound, we must consider two cases: 
\medskip 
\item{i)} The supremum is {\it finite}: $\sup_{i\geq 1} R_i < \infty$.
\smallskip 
\item{ii)} The supremum is {\it infinite}: $\sup_{i\geq 1} R_i = \infty$. 
\medskip 
{In case (i), let} $\e > 0$. There exists some $j \geq 1$ with $p_j > 0$ such that $R_j \geq (1-\e){\sup_{i \geq 1} R_i}$. We define $R = R_j$.

{In case (ii), let $R \in \R$ be an arbitrarily large value. We have $\sup_{i \geq 1} R_i = \infty$, therefore for any choice of $R$, there must exist some $j$ with $p_j > 0$ such that $R_j \geq R$.}

Now, in both cases, define the set of equivalence class indices {with a larger ratio}: 
$$ {\cal J}= \{\ell : R_{J_\ell} \geq {R}\}.$$
The probability that the {\sc MLE} chooses an equivalence class that is not a part of this set is the probability that ${\cal J}$ is empty,
$$ \pr\{\lambda \not \in {\cal J} \} = \pr\Big\{\bigcap_{\ell \in {\cal J} } N_{J_\ell} = 0\Big\} \leq \pr\{N_j^* = 0\},\nonameoldno$$
which approaches $0$ as $n \to \infty$. 
We can thus bound the probability of success from below by
\newcount\MLEprobLBfinite
\MLEprobLBfinite=\eqcount
$$
\eqalign{
\pr\{\mle\} &= \ex\big\{\pr\{\mle \mid F_n\}\big\} \cr 
&\geq \ex\left\{\ind{\lambda \in {\cal J} } {R_{J_\lambda} \over \sum_k N_{J_k}R_{J_k} } \right\} \cr
&\geq {{R} \over n(1+\eps)}  \ex\left\{\ind{\lambda \in {\cal J} }\bigind{\sum_k N_{J_k}R_{J_k} \leq n(1+\eps)}\right\}\cr
&\geq {{R} \over n(1+\eps) } \bigg(1 - \pr\{\lambda \not \in {\cal J} \} - \pr\Big\{\sum_k N_{J_k}R_{J_k} > n(1+\eps)\Big\} \bigg). \cr
}\nonameoldno$$
As $n\to\infty$, we have that $\pr\{\lambda \not \in {\cal J} \} \to 0$ and, again noting that $\sum_{k} N_{J_k} R_{J_k}$ is within $O(1)$ of the random variable $W_n = \sum_v R_{\deg^*(v)}$ defined in Lemma~\the\weightedsum, we also have $\pr\{\sum_k N_{J_k}R_{J_k}/n > 1+\e\} \to 0$. {Thus, in both cases (i) and (ii), the sum of terms in the parentheses approaches 1 as $n \to \infty$.}

{In case (i), we had $R \geq (1-\e)\sup_{i\geq 1} R_i$.}
Thus, since $\e$ was arbitrary,
$$ \liminf_{n\to\infty} n\pr\{\mle\} \sup_{i\geq 1} R_i,\nonameoldno$$
and we have equality in the limit.


{In case (ii),} 
$$\liminf_{n\to\infty} n\pr\{\mle\} \geq R\nonameoldno$$
%
for any arbitrarily large choice of $R$. We thus have
$$\lim_{n\to\infty} n \pr\{\mle\} = \infty = \sup_{i\geq 1} R_i,\nonameoldno$$
completing case (ii).\slug

This theorem applies to any distribution for which if there is a positive integer $i$ without any probability mass, then all integers $j\geq i$ have {$p_j = 0$} as well. Most of the important examples we consider satisfy this condition. We claimed earlier that in many cases, the probability of correctness is $c/n$ in the limit for some constant $c\geq 1$; indeed, Theorem \the\mainmainthm\ has shown that if there are no special nodes, then $c = \sup_{i\geq 1} R_i$ (when this is finite). In fact, since the only valid offspring distribution with mean 1 and $p_i/p_{i-1} = 1/i$ for all $i\geq 1$ is the $\Pos(1)$ distribution, the only case where $c=1$ is the family of Cayley trees, which we treated in Section \the\karysect. In most other cases, the {\sc MLE} does better, asymptotically speaking, than choosing uniformly at random.

Although the limit of $n\pr\{\mle\}$ may be infinite, the following lemma shows that it is always $o(n)$ {if no special integer is observed}. It will also apply to distributions containing special integers. We once again let $S_n$ denote the event that there exists a special node in a given free tree $F_n$, and let $\comp S_n$ denote the complement of this event.

\newcount\smallohlemma
\smallohlemma=\thmcount
\advance\thmcount by 1
\proclaim Lemma \the\smallohlemma. 
Let $T$ be a random free tree of size $n$ corresponding to a Galton--Watson tree with offspring distribution $p_i$. Let $\S$ be the set of special integers of this distribution. If $0 < \sigma^2 < \infty$, $\sup_{i\geq 1, i \not \in S} p_i/p_{i-1} < \infty$, then the probability of correctness of the {\sc MLE} satisfies 
$$\lim_{n\to\infty} \pr\{{\mle \cap {\comp S_n}}\} = 0.$$
Note that if there are no special integers in the distribution, this is exactly $\pr\{\mle\}$.

\proof For a conditional Galton--Watson tree of size $n$, recall the random variable $M_n = \max_{1 \leq i \leq n} \xi_i$ that we defined in Lemma~\the\maxdegree\ to describe the maximum degree. Next, we define $\kappa = \sup_{i\geq 1, i \not \in \S} p_i / p_{i-1} < \infty$.
Letting $\lambda \not \in \S$ be the class chosen by the {\sc MLE}, we note that the best ratio can be bounded by
$$ R_\lambda \leq \kappa (M_n+1) \leq 2\kappa M_n.$$
As for the sum of ratios over all nodes in the free tree, {note that given the event $\comp S_n$ and letting the event $W_n$ be as in \eqref{\the\Wneqdef}, we have}
$$ \sum_k N_{J_k} R_{J_k} = W_n + {{Dp_D} \over {p_{D-1}}} - {{(D+1)p_{D+1}} \over {p_{D}}},\nonameoldno$$
where $D$ is the degree of the root of the tree. Then, since $p_{D+1}/p_D \leq \kappa$,
$$ \sum_k N_{J_k} R_{J_k} \geq W_n - (D+1)\kappa .$$
Let $E_n$ be the event that $W_n \geq 2 \kappa \sqrt{n}$ and $D+1 \leq \sqrt{n}$. Observe that
$$\eqalign{
\pr\{{\comp E_n}\} &\le \pr\bigl\{ W_n < 2\kappa\sqrt n\bigr\} + \pr\bigl\{D+1\ge \sqrt n\bigr\} \cr
&\le 4\kappa^2 n \ex\bigl\{(1/W_n)^2\bigr\} + {\ex\{D+1\}\over \sqrt n} \cr
&= O\biggl({1\over n}\biggr) + {\sigma^2 + 2 + o(1)\over \sqrt n},\cr
}$$
where we used the fact that $\ex\bigl\{1/({W_n})^2 \mid |T| = n\bigr\} = O(1/n^2)$ and that
$\ex\{D\} = \sigma^2 + 1$ for the Kesten tree, to which the conditional
Galton--Watson tree locally converges. When $E_n$ holds, we have $\sum_k N_k R_K \ge W_n - \kappa\sqrt n \ge W_n/2$.


{By} Lemmas~\the\weightedsum\ and \the\weightedcontd{, we have} $\rev{W_n}/\gamma n \to 1$ in probability given $|T| = n$, and $\ex\big\{\gamma n / {W_n} \mid |T| = n \big\} \to 1$ as $n$ tends to infinity. The probability of correctness of the {\sc MLE} can thus be bounded by
\newcount\smallohbound
\smallohbound=\eqcount
$$\eqalign{
\pr\{\mle \cap {\comp S_n}\} &\le \pr\{ \mle\cap {\comp S_n} \cap E_n\} + \pr\{{\comp E_n}\} = \pr\{\mle \cap {\comp S_n} \cap E_n\} + o(1) \cr
&= \ex_{F_n}\big\{\pr\{\mle \mid F_n\} \one_{{\comp S_n} \cap E_n} \big\} +o(1)\cr
&= \ex\left\{{R_\lambda \over \sum_k N_{J_k} R_k} \one_{{\comp S_n}\cap C_n} \right\} + o(1) \cr
&\leq  \ex\left\{ {2\kappa M_n \over {W_n/2}} \ \Big| \ |T| = n \right\} + o(1)\cr 
&\leq 4\kappa \sqrt{\ex\left\{M_n^2 \mid |T| = n\right\} \ex\left\{1/{{W_n}}^2 \mid |T| = n\right\} }+o(1). 
}\nonameoldno$$
To bound $\ex\left\{M_n^2 \mid |T| = n\right\}$, let $A$ once again denote the event defined in ({\oldstyle{\the\eventAeq}});
we have
$$\eqalign{
\ex\left\{M_n^2 \mid |T| = n\right\} &= {\ex\left\{M_n^2 \one_A \right\} \over \pr\{A\} } \leq { n^2 \pr\{M_n \geq n^{7/8}\} + n^{7/4} \pr\{A\} \over \pr\{A\} } \cr 
&\leq \Theta\big(n^{5/2}\big) \pr\{M_n \geq n^{7/8}\} + n^{7/4}.\cr
}
$$
We proceed by applying the union bound to obtain
$$\eqalign{
\ex\left\{M_n^2 \mid |T| = n\right\} &\leq n\Theta\big(n^{5/2}\big) \sum_{i \geq n^{7/8}} p_i + n^{7/4} \cr 
&\leq \Theta\big(n^{7/2}\big) \sum_{i\geq 1} {i^2 p_i \over {n}^{7/4}} + n^{7/4} \cr 
&= \Theta\big(n^{7/4}\big),
}\nonameoldno$$
where the last equality follows from the fact that $\sigma^2 < \infty$. Substituting everything into ({\oldstyle{\the\smallohbound}}), we have
$$\pr\{\mle \cap {\comp S_n}\} = 2\kappa\sqrt{O\big(n^{7/4}\big)O\big(1/n^2\big)} =  O\bigg({1\over n^{1/8}}\bigg).\noskipslug\nonameoldno$$
\medskip

\boldlabel Distributions with special integers.
\hldest{xyz}{}{with}
We can now deal with the situation in which the distribution contains one or more special integers. It is clear that the {\sc MLE} should do no worse here than in the non-special case, since there is now the possibility of stumbling upon a node that must be the root.

\newcount\specialint
\specialint = \thmcount
\advance\thmcount by 1
\proclaim Theorem \the\specialint. Fix a random free tree of size $n$ corresponding to a Galton--Watson tree with offspring distribution $p_i$. Let $\S$ denote the set of special integers and suppose that $\S\neq \emptyset$, $0 < \sigma^2 < \infty$, and $\sup_{i\not \in \S} p_i/ p_{i-1} < \infty$. The probability of the {\sc MLE} being correct satisfies 
$$\lim_{n\to\infty} \pr\{\mle\} = \sum_{i\in\S} ip_i + o(1) .\nonameoldno$$

\proof The special integers $i\in\S$ satisfy $p_i \neq 0$ and $p_{i-1} = 0$. Recall from case (i) of Theorem~\the\mainmainthm\ that if there exists a node in the free tree with some special degree $i \in \S$, then there can only be one such node: $\sum_{i \in \S} N_i^* \leq 1$. Thus we can split $\pr\{\mle\}$ into two cases: Let $S_n$ and $\comp S_n$ be defined as in the previous lemma. Then
$$\eqalign{ 
\pr\{\mle\} &= \ex\big\{\pr\{\mle \mid F_n\}\big\} \cr
&= \ex \big\{ \pr\{\mle \mid F_n\} \one_{S_n} \big\} + \pr\{\mle \one_{\comp S_n}\} 
}\nonameoldno$$
The first term here is simply $\pr\{S_n\}$, since the {\sc MLE} satisfies $\pr\{\mle \mid F_n\}\one_{S_n} = \one_{S_n}$. 
As stated in the proof of Theorem~\the\mainmainthm, a conditional Galton--Watson tree converges {locally} to Kesten's limit tree as $n \to \infty$. Thus, the existence of a $u \in F_n$ with $\deg^*(u) \in \S$ is the event that a random conditional Galton--Watson tree has root of degree $i \in \S$, which occurs with probability $\sum_{i \in \S} ip_i + o(1)$. {Then, noting that $\pr\{\mle \cap {\comp S_n}\} = o(1)$ by Lemma~\the\smallohlemma, we have}
$$ \pr\{\mle\} = \sum_{i \in S} ip_i + o(1).\noskipslug\nonameoldno
$$

Comparing this result with Theorem~\the\mainmainthm, we see that the {\sc MLE} fares a lot better when there are special integers in the distribution. When there are no special integers, the product $n\pr\{\mle\}$ approaches $\sup_{i\geq 1} R_i$ (and in many cases this supremum is a constant), but we have now shown that the presence of special integers causes $\pr\{\mle\}$ itself to approach a nonzero constant.







\advance\sectcount by 1
\section\the\sectcount. Further Examples
\hldest{xyz}{}{sec\the\sectcount} 

We are now able to calculate the correctness of the {\sc MLE} for Galton--Watson trees with much more general offspring distributions. 
We hope that the examples below will demonstrate the simplicity of our general approach to deriving and analyzing the {\sc MLE}.
{A summary of these examples appears in Table 1}.
\medskip

\boldlabel Full binary trees. This is an example of a distribution with a special integer. In a full binary tree, a node either has two children or none, so we have $p_0 = p_2 = 1/2$ and 2 is a special integer. If there is only one node, then it is certainly the root. Otherwise, the root has graph-degree 2. As asserted in the previous section, there can only be one node in the free tree with graph-degree 2. In other words, for $n\geq 2$, we are guaranteed to be in case (i) of the MLE and we can choose the root with probability 1.
\medskip

\boldlabel Motzkin trees. These are also known as unary-binary trees, because every node can have either one or two children. Unlike a Catalan tree, a node can have one child in only one way, so these trees arise by the probability distribution $p_0 = p_1 = p_2 = 1/3$. When the tree has $n\geq 2$ nodes, the root has either degree 1 or 2, and we have
$$R_i = {ip_i \over p_{i-1}} = i\nonameoldno$$
for $i = 1,2$. The best strategy is to choose uniformly among all nodes with graph-degree 2, unless there are none, in which case we choose a leaf. By Theorem~\the\mainmainthm, we conclude that $n\pr\{\mle\}$ approaches 2 as $n$ gets large, so $\pr\{\mle\} \sim 2/n$.

\medskip

\boldlabel Planted plane trees. Also called {\it rooted ordered trees}, this is the family of trees that can be embedded in the plane in a unique way; reordering the subtrees of a given node produces a different tree even if these subtrees are visually indistinguishable. Random planted plane trees correspond to conditional Galton--Watson trees with a $\Geo(1/2)$ offspring distribution. Thus $p_i = 1/2^{i+1}$ for every $i$ and we have
$$R_i\bigg/\sum_v R_{\deg^*(v)} = {i/2 \over \sum_v \deg^*(v)/2} = {i\over 2(n-1)}.\nonameoldno$$
This is the probability that a node with degree $i$ is the root. The optimal strategy here is therefore to pick uniformly at random among the nodes of highest degree.

The maximal degree {$M_n$ of $T_n$} is a random variable, but we were able to give upper and lower bounds in Lemma~\the\maxdegree. For an upper bound, we have
$$\pr\{M_n \geq x\} \leq \big(1 + o(1)\big) n \pr\{\xi \geq x\} \sim n/2^x\nonameoldno$$
and this tends to 0 if $x = \log_2 n + \omega(1)$. (The small-omega notation $\omega(1)$ denotes a term {$a_n$ such that $a_n \to \infty$ as $n \to \infty$.}) Likewise, we can derive the lower bound
$$\pr\{M_n \leq x\} \leq \big(\beta + o(1)\big) \exp\big({-n}\pr\{\xi\geq x\}\big) \sim \beta\exp(-n/2^{x+1})\nonameoldno$$
for the constant $\beta$ given by Lemma~\the\maxdegree\ and this goes to 0 provided that $x = \log_2 n - \omega(1)$.
In other words,
$$\lim_{n\to \infty} \pr\{M_n \geq \log_2 n + \omega(1) \mid |T| = n\} = 0$$
and
$$\lim_{n\to \infty} \pr\{M_n \leq \log_2 n - \omega(1) \mid |T| = n\} = 0,$$
i.e., $M_n/\log_2 n \to 1$ in probability. This means that for a planted plane tree,
$$\pr\{\mle\} = {\ex\{M_n\} \over 2(n-1)} \sim {\log_2 n\over 2n}.\nonameoldno$$

\topinsert
$$\vcenter{\vbox{
\vskip -15pt

\centerline{\eightpoint {\smallheader Table 1.}
THE PROBABILITY OF CORRECTNESS OF THE MLE FOR SOME FAMILIES OF TREES}
}}$$
\vskip -30pt
$$\centerline{\vbox{
\eightpoint
\tabskip=.5em plus2em minus.5em
    \halign{
        \hfil#\hfil & \hfil#\hfil & \hfil$\displaystyle{#}$\hfil & \hfil$\displaystyle{#}$\hfil &
        \hfil#\hfil &
        \hfil$\displaystyle{#}$\hfil \cr
        \noalign{\hrule}
        \noalign{\medskip}
        Family & Distribution & R_i & {R_i\over \sum_{v} R_{\deg(v)}} & MLE & \pr\{\mle\} \cr
        \noalign{\medskip}
        \noalign{\hrule}
        \noalign{\medskip}
        $k$-ary & $\Bin(k,1/k)$ & {1\over k-1} & {k-i+1\over (k-1)n + 2} &Leaf& {k\over (k-1)n + 2} \cr
        \noalign{\medskip}
        Cayley & $\Pos(1)$ &  0 & 1/n & Choose uniformly & 1/n \cr 
        \noalign{\medskip}
        Full binary & $\Uni\{0,2\}$ & \cases{0,&if $i=1,3$;\cr \infty,&if $i=2$.} & \cases{0,&if $i=1,3$;\cr \infty,&if $i=2$.} & Degree $2$ & 1 \cr
        \noalign{\medskip}
        Planted plane & $\Geo(1/2)$ & 1/2 & {i\over 2(n-1)} &Maximize degree&{\ex\{M_n\}\over 2(n-1)}\sim{\log_2 n\over 2n} \cr
        \noalign{\medskip}
        Motzkin & $\Uni\{0,1,2\}$ & \cases{i,&if $i=1,2;$\cr 0,&if $i=3$.} & \cases{(i+o(1))/n,&if $i=1,2;$\cr 0,&if $i=3$.} &Degree 2&
       {} {2+o(1)\over n}
        \cr
        \noalign{\medskip}
        \noalign{\hrule}
    }
}}$$
\vskip -10pt
\endinsert

\medskip\boldlabel \llap{*}Large-tailed distributions. Assume that $R_i$ is strictly increasing as a function of $i$ and that $p_i/p_{i-1}\to 1$ as $i\to \infty$. For example, we may consider distributions with a polynomial tail
$$p_i = {\theta\over(i+1)^\alpha},$$
for $i\geq 1$ and $\alpha > 3$. The bound on $\alpha$ ensures that $\sigma^2<\infty$. Noting that $N_i^* / n \to p_{i-1}$, we obtain 
$$\sum_{i=1}^\infty {N^*_i R_i\over n} \to \sum_{i=1}^\infty p_{i-1} {ip_i \over p_{i-1}} = 1$$
in probability, and thus
$$\bigg|\pr\{\mle\mid F_n\} - {M\over n} \bigg| \leq f(M,n)\nonameoldno$$
where $f(M,n)/(M/n) \to 0$ in probability as $n\to \infty$. Thus we have, in general,
$$\pr\{\mle\} \sim {\ex\{M\}\over n}.\nonameoldno$$
For $p_i = \theta/(i+1)^\alpha$, we  that $\ex\{M\} = \Theta(n^{1/(\alpha-1)})$ and so our probability of correctness is $\Theta(n^{-(\alpha-2)/(\alpha-1)})$; varying $\alpha$ produces distributions with a whole range of correctness probabilities.

\section Acknowledgements

This paper was created entirely in social isolation during the 2020 coronavirus pandemic. We would like to thank the people who made Zoom, the people who invented the iPad, and our McGill comrades Konrad Anand, Jad Hamdan, Tyler Kastner, Gavin McCracken, Ndiam\'e Ndiaye, and Rosie Zhao for their support and valuable technical feedback. {We are also grateful to the two anonymous referees for suggesting changes that substantially improved the clarity and readibility of the paper.}

\medskip
\section References   
\frenchspacing

\font\tenit=cmsl10
\textfont\itfam=\tenit 

\bibliography{biblio}         
\bibliographystyle{plain}     


\bye